\documentclass[a07paper,english,17pt]{article}
\usepackage[T1]{fontenc}
\usepackage{babel}
\usepackage{graphicx,color}
\usepackage{babel}
\usepackage{amsfonts}
\usepackage{amssymb, amsmath, amsfonts, amsthm, mathrsfs}
\usepackage{hhline, latexsym}
\usepackage{float}
\usepackage[margin=3cm]{geometry}
\begin{document}
\author{Aymen Braghtha}
\title{Darboux integrable system with a triple point and Pseudo-Abelian integrals}
\maketitle
\begin{abstract}
In this paper we consider the degeneracies of the third type. More exact, the perturbations of the Darboux integrable foliation with a triple point, i.e. the case where three of the curves $\{P_i=0\}$ meet at one point, are considered. Assuming that this is the only non-genericity, we prove that the number of zeros of the corresponding pseudo-abelian integrals is bounded uniformly for close Darboux integrable foliations. Let $\mathcal{F}$ denote the foliation with triple point (assume it to be at the origin), and let $\mathcal{F}_{\lambda}=\{M_{\lambda}\frac{dH_{\lambda}}{H_{\lambda}}=0\}, M_{\lambda}$ is a integrating factor, be the close foliation. The main problem is that $\mathcal{F}_{\lambda}$ can have a small nest of cycles which shrinks to the origin as $\lambda\rightarrow0$. A particular case of this situation, namely $H_{\lambda}=(x-\lambda)^{\epsilon}(y-x)^{\epsilon_{+}}(y+x)^{\epsilon_{-}}\Delta$ with $\Delta$ non-vanishing at the origin (and generic in appropriate sense).
\end{abstract}
\text{Keywords}: \text{Abelian integrals, Limit cycles, Integrable systems}


\section{Introduction and Main Result}

This paper represent an amelioration of [Braghtha] in the sense to eleminate the technic condition on the perturbative one-form $\eta$. Precisely, we
consider an unfolding $\omega_{\lambda}$ of the one-form $\omega_0$, where $\lambda$ is a small parameter and $\omega_{\lambda}$ is a family of meromorphic one-forms 
\begin{equation}
\omega_{\lambda}=M_{\lambda}\frac{dH_{\lambda}}{H_{\lambda}},\quad H_{\lambda}=P_{\lambda}^{\epsilon}\prod^{k}_{i=1}P_i^{\epsilon_i},\quad M_{\lambda}=P_{\lambda}\prod^{k}_{i=1}P_i
\end{equation}
with $\epsilon, \epsilon_i>0, P_0, P_{\lambda}, P_i\in\mathbb{R}[x,y]$.\\

We assume that $P_0(0,0)=P_1(0,0)=P_2(0,0)=0$ and $P_i(0,0)\neq0, i=3,\ldots,k$. Generically, the triple point unfolds into three saddles $p_0^{\lambda}, p_1^{\lambda}$ and $p_2^{\lambda}$ correspond to the transversal intersections of level curves $P_{1}^{-1}(0)$ and $P_{\lambda}^{-1}(0)$, $P_{1}^{-1}(0)$ and $P^{-1}_2(0)$, and $P_2^{-1}(0)$ and $P_{\lambda}^{-1}(0)$. Here also appears a center $p_c^{\lambda}$ in the triangular region bounded by these levels curves.\\

Consider a polynomial perturbation $\omega_{\lambda, \kappa} = \omega_{\lambda} + \kappa \eta, \quad \kappa>0$ of the system $\omega_{\lambda}=M_{\lambda}\frac{\mathrm{d}H_{\lambda}}{H_{\lambda}}$, where 
$$\eta=R\mathrm{d}x +S\mathrm{d}y,$$ and $R, S\in\mathbb{R}[x,y]$ are a polynomials of degree $n$.
The foliation $\omega_{\lambda}=0$ has a maximal nest of cycles $\gamma(\lambda,h)\subseteq\{H_{\lambda}(x,y)=h\}, h\in(0,n(\lambda))$ filling a connected component of $\mathbb{R}^2\setminus\{P_{\lambda}\prod^{k}_{i=1}P_i=0\}$, which we denote $D_{(\lambda,h)}$, whose boundary is a polycycle $\gamma(\lambda,0)$. 

To such perturbation one can associate the pseudo-Abelian integral 

\begin{equation}
I(\lambda,h)=\int_{\gamma(\lambda,h)}\frac{\eta}{M_{\lambda}},
\end{equation}

which is the principal part of the Poincar\'{e} displacement function $$D(\kappa,\lambda,h)=\kappa h\int_{\gamma(\lambda,h)}\frac{\eta}{M_{\lambda}} +O(\kappa)$$
of the perturbation $\omega_{\lambda,\kappa}$ along $\gamma(\lambda,h)$.\\
 
Let us impose the following generic assumptions:\\
\linebreak
$\textbf{A}_1$ : $\frac{\partial P_{\lambda}}{\partial\lambda}|_{(0,0,0)}\neq0$.\\
$\textbf{A}_2$ : $P_1^{-1}(0), P_2^{-1}(0)$ and $P_0^{-1}(0)$ intersect transversally two by two at the origin which is the only triple point. The level curves $P_i^{-1}(0), i=3,\cdots,k$ intersect transversally and two by two.\\
\linebreak
\textbf{Theorem 1.} \emph{Let $I(\lambda,h)$ be the family of pseudo-Abelian integrals as defined above. Under assumptions $\textbf{A}_1$ and $\textbf{A}_2$, there exists a bound for the number of isolated zeros of $I(\lambda,h)=\int_{\gamma(\lambda,h)}\frac{\eta}{M_{\lambda}}$. The bound depends only on $n_i=\deg P_i, n=\max(\deg R,\deg S)$ and is uniform in the coefficients of the polynomials $P_{\lambda}, P_i, R, S$, the exponents $\epsilon,\epsilon_i, i=1,\ldots,k$ and the parameter $\lambda$.}
\subsection{Bobie\'nski result}
We recall also a similar result of Bobi\'enski [1] which he prove the existence of a local upper bound for the number of zeros of pseudo-abelian integrals. He consider a one-parameter unfolding of the singular codimension one case. The difference relies in the fact that in this work the Darboux first integrals is more general ($\epsilon_1\neq\epsilon_2$) and the proof of our main result is purely geometric: we use the blow-up in families. This approach gives directly uniform validity of our study of the pseudo-abelian integrals.

\section{Darboux integrable foliation}
Let us a stablish a local normal form near the triple point $(0,0,0)$ for the unfolding of the degenerate polycyle  $H_{0}$.\\
\linebreak
\textbf{Proposition 1.} \emph{Under above assumptions $\textbf{A}_1, \textbf{A}_2$. There exists a local analytic coordinate system $(x,y,\lambda)$ at $(0,0,0)$ such that $H_{\lambda}$ takes the form
\begin{equation}
H_{\lambda}=(x-\lambda)^{\epsilon}(y-x)^{\epsilon_{+}}(y+x)^{\epsilon_{-}}\Delta,\quad\lambda>0
\end{equation}
where $\Delta$ is an analytic unity function $\Delta(0,0,0)\neq0$.}\\
\linebreak
Let $\mathcal{F}_1, \mathcal{F}_2$ are two foliations of dimension two in space of the total complex space $\mathbb{C}^{3}$ with coordinates $(x,y,\lambda)$ such that 
$$\mathcal{F}_1:\{H(x,y,\lambda)=P_{\lambda}^{\epsilon}\prod^k_{i=1}P_i^{\epsilon_i}=h\},\quad
\mathcal{F}_2:\{\lambda=\text{constant}\}.$$

Consider the Darboux foliation $\mathcal{F}:=\{\mathcal{F}_1,\mathcal{F}_2\}$ of dimension one in $\mathbb{C}^3$ with coordinates $(x,y,\lambda)$ which is given by the intersection of the leaves of $\mathcal{F}_1$ and $\mathcal{F}_2$. This foliation has a non-elementary singular point at the origin $(0,0,0)$. This foliation has a complicated singularity at the origin $(0,0,0)$. To reduice this singularity , we perform a directional blowing-up in the family $\omega_{\lambda}$. Note that blow-up in a family was introduced in [5], see also [6].
\subsection{Directional Blow-up}
The blow-up of $\mathbb{C}^{3}$ at the origin is defined as the incidence three dimensional manifold $W=\{(p,q)\in\mathbb{CP}^{2}\times\mathbb{C}^{3}: q\in p\}$. The blow down $\sigma:W\rightarrow\mathbb{C}^{3}$ is just the restriction to $W$ of the projection $\mathbb{CP}^2\times\mathbb{C}^{3}$. The inverse map $\sigma^{-1}:\mathbb{C}^{3}\rightarrow W$ is called blow-up and $\sigma^{-1}(0)=\mathbb{CP}^2$ is called exceptional divisor. The projective space $\mathbb{CP}^{2}$ is covered by three canonical charts: $W_1 = \{x\neq0\}$ with coordinates $(v_{1}, w_{1})$, $W_{2} = \{y \neq 0\}$ with coordinates $(u_{2}, w_{2})$ and $W_{3} =\{\lambda\neq0\}$ with coordinates $(u_3,v_3)$. $W_1, W_2$ and $W_3$ define canonical charts on $W$, with coordinates $(u_1, v_1, w_1)$,\\
$(u_2, v_2, w_2)$ and $(u_3, v_3, w_3)$ respectively. The blow-up $\sigma^{-1}$ is written as:
\begin{align}
&\sigma^{-1}_{1}=\sigma^{-1}|_{W_{1}}: x = u_1  \qquad  y = u_1 w_1  \qquad  \lambda = w_1 u_1\\
&\sigma^{-1}_{2}=\sigma^{-1}|_{W_{2}}: x = u_2 v_2  \qquad  y = v_2  \qquad  \lambda = w_2 v_2\\
&\sigma^{-1}_{3}=\sigma^{-1}|_{W_{3}}: x = u_3 w_3 \qquad y = v_3 w_3 \qquad \lambda = w_3
\end{align}


Let $\sigma^{-1}\mathcal{F}$ be the lift of the foliation $\mathcal{F}$ to the complement of the exceptional divisor $\mathbb{CP}^2$. This foliation is regular outside of the preimage of the hypersurface $\{H_{\lambda}=0, \lambda=0\}$. \\
\linebreak
\textbf{Prposition 2.} \emph{Let $a=\epsilon+\epsilon_-+\epsilon_+$
\begin{enumerate}
\item The foliation $\sigma^{-1}\mathcal{F}$ extends in a unique way to a holomorphic singular foliation $\sigma^{\ast}\mathcal{F}$ on $W$ which we call the blow-up of the original codimension two foliation $\mathcal{F}$ by the map $\sigma$.\\
\item Let $\sigma^{\ast}_1\mathcal{F}$ be the restriction of the blown-up foliation $\sigma^{\ast}\mathcal{F}$ to the chart $W_1$. The singularities of $\sigma^{\ast}_1\mathcal{F}$ are located at the points $p_+=(0,1,0), p_-=(0,-1,0), q_+=(0,1,1)$ and $q_-=(0,-1,1)$.  All these singular points are linearisable saddles, with eigenvalues $\mu_+=(\epsilon_+,-a,-\epsilon_-), \mu_-=(-\epsilon_-,a,\epsilon_-), \nu_+=(0,-\epsilon,\epsilon_+)$ and $\nu_{-}=(0,-\epsilon,\epsilon_-)$ respectively.
\end{enumerate}}
\begin{proof} We prove the second item of Proposition. Since $\sigma:W\rightarrow\mathbb{C}^{3}$ is a biholomorphism outside $\mathbb{CP}^2$, all singularities of $\sigma^{\ast}_1\widetilde{\mathcal{F}}$ on $W_1\setminus\{u_1=0\}$ correspond to singularities of $\mathcal{F}$.

Thus, it suffices to compute the singularities of $\sigma^{\ast}_1\mathcal{F}$ on the exceptional divisor $\{u_1=0\}$. On the exceptional divisor, the foliation is given by the levels of  
$$
G:=\frac{(\pi\circ\sigma_1)^{a}}{H\circ\sigma_1}=w_1^{a}(1-w_1)^{-\epsilon}(v_1-1)^{-\epsilon_+}(v_1+1)^{-\epsilon_-}\tilde{\Delta}^{-1},
$$
where $\tilde{\Delta}$ is unit of the form $\tilde{\Delta}=c+u_1f, f$ is a holomorphic function. 
 Let us compute the eigenvalues at $p_+, p_-, q_+$ and $q_-$. Near the exceptional divisor $\{u_1=0\}$, the foliation $\sigma_1^{\ast}\mathcal{F}$ is given by 
\begin{align*}
&H\circ\sigma_1(u_1,v_1,w_1)=u_1^{a}(1-w_1)^{\epsilon}(v_1-1)^{\epsilon_+}(v_1+1)^{\epsilon_-}\tilde{\Delta}=h,\\
&\pi\circ\sigma_1(u_1,v_1,w_1)=u_1w_1=\lambda.
\end{align*}
  
Near $p_{+}$ and after the respective changes of variable $v=(v_1-1)(v_1+1)^{\frac{\epsilon_{-}}{\epsilon_{+}}}\Delta$, the blown-up foliation $\sigma^{\ast}_1\mathcal{F}$ is given by the two first integrals $u_1^{a}v_{+}^{\epsilon_{+}}=h$ and $u_1w_1=s$. Then, near this point the vector field generating the foliation $\sigma^{\ast}_1\mathcal{F}$ is given by
$$
X(u_1,v,w_1)=\mu_{1}^{+}u_1\frac{\partial}{\partial u_1}+\mu_2^{+}v\frac{\partial}{\partial v}+\mu_3^{+}w_1\frac{\partial}{\partial w_1},
$$
where the vector $\mu_{\pm}=(\mu^{\pm}_1,\mu^{\pm}_2,\mu^{\pm}_3)$ satisfies the following equations
$$
<(\mu_1^{\pm},\mu_2^{\pm},\mu^{\pm}_3),(a,\epsilon_{\pm},0)>=0,\quad <(\mu^{\pm}_1,\mu_2^{\pm},\mu_3^{\pm}),(1,0,1)>=0.
$$
Here $<,>$ is the usual scalar product on $\mathbb{C}^3$. By simple calculations, we obtain 
$$
X_{\pm}(u_1,v_{\pm},w_1)=\pm\epsilon_{\pm}u_1\frac{\partial}{\partial u_1}\mp av_{\pm}\frac{\partial}{\partial v_{\pm}}\mp\epsilon_{\pm}w_1\frac{\partial}{\partial w_1}.
$$
Similar computation shows that there are local coordinates near $q_{\pm}$ in which the vector field generating the foliation is given by
$$
X_{\pm}(u_1,v_{\pm},w_{\pm})=-\epsilon w_{\pm}\frac{\partial}{v_{\pm}}+\epsilon_{\pm} w_{\pm}\frac{\partial}{\partial w_{\pm}}.
$$
\end{proof} 
\section{Proof of the theorem}
Let $t=\frac{\lambda^{a}}{h}$. The blown-up foliation $\sigma^{\ast}_1\mathcal{F}$ is given by the two first integrals 
\begin{align*} 
&G=w_1^{a}(1-w_1)^{-\epsilon}(v_1-1)^{-\epsilon_+}(v_1+1)^{-\epsilon_-}\tilde{\Delta}^{-1}=t,\\
&F=u_1w_1=\lambda.
\end{align*}  
Consider the two-dimensional square $Q\subset\mathbb{CP}^2$ with vertices $p_+,p_-,q_+$ and $q_-$. All levels curves $\{G=t\}$ inside $Q$ correspond to values of $t\in[0,+\infty]$. We consider the family of hyperbolic polycycles $\delta^t$, i.e. at each intersection of two consecutive curve we have a saddle point,.
$$
\delta^{t}=\left(\sigma^{-1}_1\left(\gamma(0,0)\setminus(0,0,0)\right)\cup\left(Q\cap\{G=t\}\right)\right)^{\mathbb{R}}, t\in[0,+\infty].
$$
where $\left(\ldots\right)^{\mathbb{R}}$ denotes the real part of a complex analytic set. All polycycle $\delta^t$ satisfies the generecity assumptions from [2].


Let $\delta(\lambda,t)=\sigma^{-1}(\gamma(\lambda,h))\subset W$ be the pull-back of the cycle $\gamma(\lambda,h)$ by the blowing-up map. Let $\delta^t$ be the polycycle corresponding to the cycle $\delta(\lambda,t)$.
Let
$$
J(\lambda,t)=\int_{\delta(\lambda,t)}\sigma^{\ast}_1\frac{\eta}{M_{\lambda}}.
$$
The integral $J(\lambda,t)$ can be analytically continued to the universal cover of $\mathbb{C}^{2}\setminus\{\lambda t=0\}$.
\subsection{Variation operator}
Given any multivalued function $F$ defined in a neighborhood of the origin in $\mathbb{C}$ i.e. a holomorphic function defined on the universal covering $\widetilde{\mathbb{C}^{\ast}}$ of $\mathbb{C}^{\ast}$. We define the rescaled monodromy as
$$
\mathcal{M}on_{(t,\alpha)}F(t)=F(te^{i\pi\alpha}).
$$
The variation is given as the difference between the counterclockwise and clockwise continuation 
\begin{align*}
\mathcal{V}ar_{(t,\alpha}F(t)&=\mathcal{M}on_{(t,\alpha)}J(t)-\mathcal{M}on_{(t,-\alpha)}J(t)\\
&=F(te^{i\pi\alpha})-F(te^{-i\pi\alpha}). 
\end{align*}
Now, let $G$ be a multivalued function in two variables $\lambda$ and $t$ defined in universal covering $\widetilde{\mathbb{C}^2\setminus\{st=0\}}$ of $\mathbb{C}^2\setminus\{st=0\}$. We define the mixed variation as 
\begin{align*}
&\mathcal{V}ar_{(\beta,\lambda)}\circ\mathcal{V}ar_{\alpha,t}G(\lambda,t)=\mathcal{V}ar_{(\beta,\lambda)}(G(\lambda,te^{i\pi\alpha})-G(\lambda,te^{-i\pi\alpha}))=\\
&G(\lambda e^{i\pi\beta},te^{i\pi\alpha})-G(\lambda e^{-i\pi\beta},te^{i\pi\alpha})-G(\lambda e^{i\pi\beta},te^{-i\pi\alpha})+G(\lambda e^{-i\pi\beta},te^{-i\pi\alpha}).
\end{align*}
\textbf{Lemma 1} \emph{The variations $\mathcal{V}ar_{(\lambda,\beta)}$ and $\mathcal{V}ar_{(t,\alpha)}$ commute
$$
\mathcal{V}ar_{(\lambda,\beta)}\circ\mathcal{V}ar_{(t,\alpha)}=\mathcal{V}ar_{(t,\alpha)}\circ\mathcal{V}ar_{(\lambda,\beta)}.
$$ }
\begin{proof}
The proof is a direct consequence of monodromy theorem.
\end{proof}
\subsection{The analytic proporties of the integral $J$ near the polycycle $\delta^t$}
Using the partition of unity of the blown-up space we can decompose our cycle of integration $\delta(\lambda,t)$ in a relative cycles and we can check that any relative cycle can be chosen as a lift of a base path (for more details see [3]). Then, we have\\
\linebreak
\textbf{Proposition 3 .} \emph{The integral $J(\lambda,t)$ satisfies the rescaled iterated variation equations 
\begin{align}
&\mathcal{V}ar_{(t,\beta_1)}\circ\ldots\circ\mathcal{V}ar_{(t,\beta_k)}J(\lambda,t)=0,\\
&\mathcal{V}ar_{(\lambda,1)}\circ\mathcal{V}ar_{(\lambda,1)}J(\lambda,t)=0,
\end{align} where $\beta_i$ are analytic functions in $\epsilon,\epsilon_+,\epsilon_-,\ldots,\epsilon_k$. }
\subsection{Increment of argument of $J$}
The estimation of a local bound of the number of isolated zeros of the pseudo-abelian integral $I(\lambda,h)=\int_{\gamma(\lambda,h)}\frac{\eta}{M_{\lambda}}$ is analoguous to estimate a local $\lambda$-uniform bound of isolated zeros of the integral $J(\lambda,t)=\int_{\delta(\lambda,t)}\sigma_1^{\ast}\frac{\eta}{M_{\lambda}}$. As consequence of the equations (7) and (8) the integral $J(s,t)$ has the following expansion
\begin{equation}
J(\lambda,t)=J_1(\lambda,t)+J_2(\lambda,t)\log\lambda,
\end{equation}
where 
\begin{align}
&J_2(\lambda,t)=\mathcal{V}ar_{(\lambda,1)}J(\lambda,t)=\int_{\text{eight loop}}\sigma^{\ast}_1\frac{\eta}{M_{\lambda}},\\
&\mathcal{V}ar_{(t,\beta_1)}\circ\ldots\circ\mathcal{V}ar_{(t,\beta_k)}J_i(\lambda,t)=0,\quad i=1,2,\\
&\mathcal{V}ar_{(\lambda,1)}J_i(\lambda,t)=0, i=1,2.
\end{align}
\textbf{Theorem 2.} \emph{Let $\epsilon>0$ be sufficiently small. Then, for all  $|\lambda|<\epsilon$ the number of zeros $$\#\{t\in[0,+\infty]: J(\lambda,t)=0\}$$ is uniformly bounded with respect to $\lambda$.}\\
\linebreak
The integral $J(\lambda,t)$ has analytic prolongation to the complex argument $t$. This is a multivalued function with unique ramification point $t=0$.\\

Let $\partial \Omega$ be the boundary of a domain complex $\Omega$ which consists of a big circular arc $C_{R_1}=\{|t|=R_1, |\arg t|\leq\alpha\pi\}$, a two segments $C^{\pm}=\{r_1\leq|t|\leq R_1, |\arg t|=\pm\alpha\pi\}$ and the small circular arc $C_{r_1}=\{|t|=r_1, |\arg t|\leq\alpha\pi\}$.

To count the number of zeros of the function $J(s,t)$, we estimate the increment of argument of the function $J(\lambda,t)$ we apply the argument principle which says that
$$
\#Z(J_{\Omega})\leq\frac{1}{2\pi}\Delta\arg(J_{\partial \Omega})=\frac{1}{2\pi}(\Delta\arg(J_{C_{R_1}})+\Delta\arg(J_{C_{r_1}})+\Delta\arg(J_{C^{\pm}}))
$$
The increment of arguments $\Delta\arg(J_{C_{R_1}})$ and $\Delta\arg(J_{C^{\pm}})$ are locally uniformly  bounded, for more details see [3]. The problem consist to estimate the increment of argument of the integral $J(s,t)$ along the small circular arc $C_{r_1}$.

Let 
\begin{align*}
\Lambda(\lambda;n_0,n_1,\cdots,n_k;n)=&\{M_{\lambda}\frac{dH_{\lambda}}{H_{\lambda}}+\kappa(Rdx+Rdy): H_{\lambda}=P_{\lambda}\prod^k_{i=1}P_i^{\epsilon_i},\\
& \deg P_{\lambda}\leq n_0, \deg P_i\leq n_i, \deg(R,S)\leq n\}.
\end{align*}
be the parameters space. Consider the following functional space $\mathcal{P}$
$$
\mathcal{P}(v,V;\alpha_1,\ldots,\alpha_k;\lambda):=\{\sum\sum c_{jln}(\lambda)t^{\alpha_jn}\log^{n}t, c_{jln}\in\mathbb{C}, v\leq\alpha_{j}n\leq V,0\leq l\leq k\}.
$$
\subsubsection{Finite order approximation to $J(\lambda,t)$}
Our goal is to obtain an asymptotic expansion for $J(\lambda,t)$. This will allow us to prove the existence of an upper bound for the increment of the argument $\Delta\arg(J_{C_{r_1}})$.
The problem is proposed in the fact that the existence of the term $\log(\lambda)\rightarrow\infty$ as $\lambda\rightarrow0$ in the expression of the function $J(\lambda,t)$.
To resolve the problem we make a weighted quasi-homogeneous blowing-up. Let 
$$
c=(\log\lambda)^{-1},\quad a=J_1(\lambda,t),\quad b=J_2(\lambda,t).
$$
Then, we have 
$$
J(\lambda,t)=J_1(\lambda,t)+\log(\lambda)J_2(\lambda,t)=c^{-1}(ca+b)=c^{-1}\psi(a,b,c).
$$
As $c^{-1}=\log\lambda\in\mathbb{R}$, we have
$$
\arg(J(\lambda,t))=\arg(c^{-1}\psi(a,b,c))=\arg(\psi(a,b,c)).
$$
\textbf{Proposition 4.} \emph{There exist a some meromorphic function $\phi$ such that  $$\Delta\arg(\psi(a,b,c))=\Delta\arg(\phi).$$ }
\begin{proof} 
To estimate the increment of argument $\arg(q(a,b,c))$ of polynomial $q$ we make a quasi-homogeneous blowing-up of weight $(\frac{1}{2},1,\frac{1}{2})$ with the center $C_1=\{a=b=c=0\}$.
The explicit formulae of the blowing-up $\pi_1^{-1}$ in the affine charts $\tau_1=\{a\ne0\}, \tau_2=\{b\ne0\}$ and $\tau_3=\{c\ne0\}$ is written respectively as
\begin{align*}
&\pi_{1}^{-1}|_{\tau_1}=\pi_{11}^{-1}: a=\sqrt{a_1},\quad b=b_1a_1,\quad c=c_1\sqrt{a_1},\\
&\pi_{1}^{-1}|_{\tau_2}=\pi_{12}^{-1}: a=a_2\sqrt{b_2},\quad b=b_2,\quad c=c_2\sqrt{b_2},\\
&\pi_{1}^{-1}|_{\tau_3}=\pi_{13}^{-1}: a=a_3\sqrt{c_3},\quad b=b_3c_3,\quad c=\sqrt{c_3}.
\end{align*}
The total transform $\pi^{\ast}_{1}q(a,b,c)$ of $q(a,b,c)$ is given, in different charts, by
\begin{align*}
&\pi^{\ast}_{11}(ca+b)=a_1(c_1+b_1)=a_1P_1(a_1,b_1,c_1),\\
&\pi^{\ast}_{12}(ca+b)=b_2(a_2c_2+1)=b_2P_2(a_2,b_2,c_2),\\
&\pi^{\ast}_{13}(ca+b)=c_3(a_3+b_3)=c_3P_3(a_3,b_3,c_3).
\end{align*}
where $\{a_1=0\},\{b_2=0\}$ and $\{c_3=0\}$ are local equations of the exceptional divisor and $\{P_1=0\},\{P_2=0\}$ and $\{P_3=0\}$  local equations of strict transform of $\{ca+b=0\}$.

We observe that the exceptional divisor $\{a_1=0\}$( resp $\{c_3=0\}$) has not a normal crossing with the strict transform $\{P_1=0\}$ (resp $\{P_3=0\}$). Conversely in the chart $\tau_2$, the exceptional divisor $\{b_2=0\}$ has a normal crossing with the strict transform $\{P_2=0\}$. To resolve this problem we make a second blowing-up $\pi_2^{-1}$ with center $C_2$ such that 
\begin{enumerate}
\item in the chart $\tau_1$ with coordinates $(a_1,b_1,c_1)$, the center $C_2$ is given by  $C_2=\{c_1=b_1=0\}$,
\item in the chart $\tau_2$ the blowing-up $\pi_2$ is a biholomorphism,
\item in the chart $\tau_3$ with coordinates $(a_3,b_3,c_3)$, the center $C_2$ is given by  $C_2=\{a_3=b_3=0\}$,
\end{enumerate}
and 
\begin{enumerate}
\item in the chart $\tau_1$, we have $\pi_2^{-1}(C_2)$ is covred by two coordinates charts $V_{b_1}$ and $V_{c_1}$ with coordinate $(\tilde{a}_1,\tilde{b}_1,\tilde{c}_1)$ such that the blowing-up $\pi_2^{-1}$ is given  in the charts $V_{b_1}$ and $V_{c_1}$ respectively, by $a_1=\tilde{a}_1, b_1=\tilde{b}_1, c_1=\tilde{b}_1\tilde{c}_1$ and $a_1=\tilde{a}_1, b_1=\tilde{b}_1\tilde{c}_1, c_1=\tilde{c}_1$,
\item in the chart $\tau_3$, we have $\pi_2^{-1}(C_2)$ is covred by two coordinates charts $V_{a_3}$ and $V_{b_3}$ with coordinates $(\tilde{a}_3,\tilde{b}_3,\tilde{c}_3)$ such that the blowing-up $\pi_2^{-1}$ is given in $V_{a_3}$ and $V_{b_3}$ respectively, by $a_3=\tilde{a}_3, b_3=\tilde{a}_3\tilde{b}_3, c_3=\tilde{c}_3$ and $a_3=\tilde{a}_3\tilde{b}_3, b_3=\tilde{a}_3, c_3=\tilde{c}_3$.
\end{enumerate}
Then, we have
\begin{enumerate}
\item in the chart $V_{b_1}$,  
$$
\pi_2^{\ast}\circ\pi_1^{\ast}\psi(a,b,c)\overset{0}{\approx}\tilde{a}_1\tilde{b}_1=J_2(\lambda,t),
$$
\item in the chart $V_{c_1}$, 
$$
\pi_2^{\ast}\circ\pi_1^{\ast}\psi(a,b,c)\overset{0}{\approx}\tilde{a}_1\tilde{c}_1=\frac{J_1(\lambda,t)}{\log\lambda},
$$
\item in the chart $V_{a_3}$, 
$$
\pi_2^{\ast}\circ\pi_1^{\ast}\psi(a,b,c)\overset{0}{\approx}\tilde{a}_3\tilde{c}_3=\frac{J_1(\lambda,t)}{\log\lambda},
$$
\item in the chart $V_{b_3}$, 
$$
\pi_2^{\ast}\circ\pi_1^{\ast}\psi(a,b,c)\overset{0}{\approx}\tilde{b}_3\tilde{c}_3=\frac{J_1(\lambda,t)}{\log\lambda},
$$
\item in the chart $\tau_2$, the exceptinal divisor $\{b_2=0\}$ has a normal crossings with the strict transform $P_2=0$ and the blowing-up $\pi_2^{-1}$ is a biholomophism, so the function 
$$
\frac{J_1(\lambda,t)}{\log\lambda}+J_2(\lambda,t)=\frac{J(\lambda,t)}{\log\lambda}
$$
is meromorphic function.
\end{enumerate}
Then, we conclude that $\phi\in\{\frac{J(\lambda,t)}{\log\lambda},\frac{J_1(\lambda,t)}{\log\lambda},J_2(\lambda,t)\}$.
\end{proof}
\subsubsection{Proof of the theorem 2}
Using lemma 4.8 of [2] for $\phi\in\{\frac{J_1(\lambda,t)}{\log\lambda}, J_2(\lambda,t)\}$, there exist a $\tilde{\phi}\in\mathcal{P}(...)$ such that $|\phi(\lambda,t)-\tilde{\phi}(\lambda,t)|\leq N$.
For $\phi=\frac{J(\lambda,t)}{\log\lambda}$, we can conclude using Gabrielov's theorem [4].

\textbf{Burgundy University, Burgundy Institue of Mathematics,\\
U.M.R. 5584 du C.N.R.S., B.P. 47870, 21078 Dijon \\
Cedex - France.}\\
\textbf{E-mail adress:} aymenbraghtha@yahoo.fr



\begin{thebibliography}{9}
\bibitem{} Bobie\'{n}ski, Marcin Pseudo-Abelian integrals along Darboux cycles\textit{ a codimension one case. J. Differential Equations 246 (2009), no. 3, 1264-1273.}
\bibitem{} Bobie\'{n}ski, Marcin; Marde\v{s}i\'{c}, Pavao \textit{Pseudo-Abelian integrals along Darboux cycles. Proc. Lond. Math. Soc. (3) 97 (2008), no. 3, 669-688.}
\bibitem{} Braghtha Aymen (2013) \textit{Les z\'eros des int\'egrales pseudo-abeliennes: cas non g\'en\'erique. PhD thesis. Universit\'e de Bourgogne: France.}
\bibitem{} Gabri\`{e}lov, A. M  \textit{Projections of semianalytic sets. (Russian) Funktsional. Anal. i Prilo\v{z}en. 2 1968 no. 4, 18-30}
\end{thebibliography}
\end{document}